\documentclass[12pt]{article}
\usepackage{amsmath}
\usepackage{amstext}
\usepackage{amsfonts}
\usepackage{amsthm}
\usepackage{amscd}
\numberwithin{equation}{section}

\swapnumbers
\theoremstyle{plain}
\newtheorem{thm}{Theorem}[section]
\newtheorem{lem}[thm]{Lemma}
\newtheorem{prop}[thm]{Proposition}

\theoremstyle{remark}
\newtheorem{rem}[thm]{Remark}

\def\mod0{{\mathcal SU}_C(2,{\mathcal O})}

\def\notin{\in \hspace{-4mm}/\ }

\def\cO{{\mathcal O}}

\def\map#1{\ \smash{\mathop{\longrightarrow}\limits^{#1}}\ }

\def\sing{\mathrm{Sing} \: \Theta}

\def\LL{\mathbb{L}}
\def\L{\mathcal{L}}

\def\PP{\mathbb{P}}

\def\TT{\mathbb{T}}

\def\EE{\mathcal{E}}

\def\CC{\mathbb{C}}

\def\D{\mathcal{D}}

\def\VV{\mathcal{V}}
\def\F{\mathbf{F}}
\def\m{\mathfrak{m}}
\def\pp{\mathfrak{p}}
\def\B{\mathcal{B}}

\def\supp{\mathrm{supp} \:}
\def\coker{\mathrm{coker} \:}
\def\deg{\mathrm{deg} \:}

\def\sym{\mathrm{Sym}}
\def\ker{\mathrm{ker} \:}
\def\pic{\mathrm{Pic}}
\def\mult{\mathrm{mult} \:}

\def\dim{\mathrm{dim} \:}
\def\det{\mathrm{det} \:}
\def\rk{\mathrm{rk} \:}

\def\im{\mathrm{im} \:}
\def\Gr{\mathrm{Gr}}
\def\Spec{\mathrm{Spec}}

\def\Rat{\mathrm{Rat}}
\def\Prin{\mathrm{Prin}}
\def\Hess{\mathrm{Hess}}
\def\mod{\mathrm{mod}}
\def\exp{\mathrm{exp}}
\def\Bs{\mathrm{Bs}}

\def\lra{\longrightarrow}

\def\ra{\rightarrow}

\def\lms{\longmapsto}

\def\cc{\mathbb{C}}

\begin{document}

\title{On cubics and quartics through a canonical curve}
\author{Christian Pauly}
\maketitle

\begin{abstract}
We construct families of quartic and cubic hypersurfaces through a
canonical curve, which are parametrized by an open subset in a
Grassmannian and a Flag variety respectively. Using
G. Kempf's cohomological obstruction theory, we show that
these families cut out the
canonical curve and that the quartics are birational (via a
blowing-up of a linear subspace) to quadric bundles over the
projective plane, whose Steinerian curve equals the canonical curve.
\end{abstract}

\section{Introduction}
Let $C$ be a smooth nonhyperelliptic curve of genus $g \geq 4$, which
we consider as an embedded curve $\iota_\omega :
C \hookrightarrow \PP^{g-1}$ by its canonical linear series $|\omega|$.
Let $I = \bigoplus_{n \geq 2} I(n)$ be
the graded ideal of the canonical curve. It was classically known
(Noether-Enriques-Petri theorem, see e.g. \cite{acgh} p. 124) that the
ideal $I$ is generated by its elements of degree $2$, unless $C$ is trigonal
or a plane quintic.

\bigskip

It was also classically known how to construct some distinguished
quadrics in $I(2)$. We consider a double point of the
theta divisor $\Theta \subset \pic^{g-1}(C)$, which corresponds by
Riemann's singularity theorem to a
degree $g-1$ line bundle $L$ satisfying $\dim |L| = \dim
|\omega L^{-1}| = 1$ and we observe that the morphism
$\iota_L \times \iota_{\omega L^{-1}} : C \lra C' \subset |L|^* \times
|\omega L^{-1}|^* = \PP^1 \times \PP^1$ (here $C'$ denotes the image curve)
followed by the Segre embedding
into $\PP^3$ factorizes through the canonical space
$|\omega|^*$, i.e.,
$$
\begin{array}{ccc}
C & \hookrightarrow & |\omega|^* \\
\downarrow & &  \downarrow^{\pi}  \\
\PP^1 \times \PP^1 & \hookrightarrow & \PP^3,
\end{array}
$$
where $\pi$ is projection from a $(g-5)$-dimensional vertex $\PP V^\perp$
in $|\omega|^*$.
We then define the quadric $Q_L := \pi^{-1}(\PP^1 \times \PP^1)$,
which is a rank
$\leq 4$ quadric in $I(2)$ and coincides with the projectivized tangent
cone at the double point $[L] \in \Theta$ under the identification
of $H^0(C,\omega)^*$ with the tangent space $T_{[L]} \pic^{g-1}(C)$.
The main result, due to M. Green \cite{green}, asserts that the set of
quadrics $\{ Q_L \}$, when $L$ varies over the double points of $\Theta$,
linearly spans $I(2)$. From this result one infers a constructive Torelli
theorem by intersecting all quadrics $Q_L$ --- at least for $C$ general
enough.

\bigskip

The geometry of the theta divisor $\Theta$ at a double point $[L]$
can also be exploited to produce higher degree elements in the ideal
$I$ as follows: we expand in a suitable set of coordinates a local
equation $\theta$ of $\Theta$ near $[L]$ as $\theta = \theta_2 +
\theta_3 + \ldots$, where $\theta_i$ are homogeneous forms of degree $i$.
Having seen that $Q_L = \mathrm{Zeros}(\theta_2)$, we denote by $S_L$ the
cubic $\mathrm{Zeros}(\theta_3) \subset |\omega|^*$,
the osculating cone of $\Theta$ at $[L]$. The cubic $S_L$ has many nice
geometric properties: under the blowing-up of the vertex $\PP V^\perp
\subset S_L$, the cubic $S_L$ is transformed into a quadric bundle
$\tilde{S}_L$ over $\PP^1 \times \PP^1$ and it was shown by
G. Kempf and F.-O. Schreyer \cite{ks}  that the Hessian and
Steinerian curves of $\tilde{S}_L$ are $C' \subset \PP^1 \times \PP^1$ and
$C \subset |\omega|^*$ respectively,
which gives another proof of Torelli's theorem.

\bigskip

In this paper we construct and study distinguished cubics and
quartics in the ideal $I$ by adapting the methods of \cite{ks} to
rank-$2$ vector bundles over $C$. Our construction basically goes as
follows (section 2): we consider a general $3$-plane $W \subset
H^0(C, \omega)$ and define the rank-$2$ vector bundle $E_W$ as
the dual of the kernel of the evaluation map in $\omega$ of
sections of $W$. The bundle $E_W$ is stable and admits a theta divisor
$D(E_W)$ in the Jacobian $JC$. Since $D(E_W)$ contains the origin
$\cO \in JC$ with multiplicity $4$, the projectivized tangent cone to
$D(E_W)$ at $\cO$ is a quartic hypersurface in $\PP T_\cO JC = |\omega|^*$,
denoted by $F_W$ and which contains the canonical curve. We
therefore obtain a rational map from the Grassmannian $\Gr(3,H^0(\omega))$
to the ideal of quartics $|I(4)|$
\begin{equation}\label{fmap4}
\F_4 \  :  \  \Gr(3,H^0(\omega)) \dashrightarrow |I(4)|, \qquad W \mapsto F_W.
\end{equation}

Our main tool to study the tangent cones $F_W$ is G. Kempf's
cohomological obstruction
theory \cite{kempf1},\cite{kempf2},\cite{ks} which in our set-up leads
to a simple criterion (Proposition \ref{simpl}) 
for $b \in \PP T_\cO JC = |\omega|^*$ to belong to $F_W$. 
We deduce in particular from this criterion
that the cubic polar $P_x(F_W)$ of $F_W$ with respect to a point
$x \in W^\perp$ also contains the canonical curve. Here
$W^\perp$ denotes the annihilator of $W \subset H^0(\omega)$.
We therefore obtain a
rational map from the flag variety $\mathrm{Fl}(3,g-1,H^0(\omega))$
parametrizing pairs $(W,x)$ to the ideal of cubics $|I(3)|$
\begin{equation}\label{fmap3}
\F_3 \  :  \  \mathrm{Fl}(3,g-1,H^0(\omega)) \dashrightarrow |I(3)|,
\qquad (W,x) \mapsto P_x(F_W).
\end{equation}

\bigskip

Our two main results can be stated as follows.

\bigskip

{\bf (1)}  Like the cubic osculating cones $S_L$, the quartic
tangent cones $F_W$ transform under the blowing-up of the vertex
$\PP W^\perp \subset F_W$ into a quadric bundle $\tilde{F}_W \ra
\PP W^* = \PP^2$. Their Hessian and Steinerian curves are the
plane curve $\Gamma$, image under the projection with center $\PP
W^\perp$, $\pi: C \ra \Gamma \subset \PP W^*$, and the canonical
curve $C \subset |\omega|^*$ (Theorem \ref{mainthm}). This surprising
analogy with the osculating cones $S_L$ remains however
unexplained.

\bigskip

{\bf (2)} Let us denote by $|F_4| \subset |I(4)|$ and $|F_3| \subset
|I(3)|$ the linear subsystems spanned by the quartics $F_W$ and the cubics
$P_x(F_W)$ respectively. Then we show
(Theorem \ref{mainthm2}) that both base loci of $|F_4|$ and $|F_3|$ coincide
with $C \subset |\omega|^*$,i.e., the quartics $F_W$ (resp. the
cubics $P_x(F_W)$) cut out the canonical curve.

\bigskip

The starting point of our investigations was the question asked by
B. van Geemen and G. van der Geer (\cite{vgvg} page 629) about ``these
mysterious quartics'' which arise as tangent cones to $2\theta$-divisors
in the Jacobian having multiplicity $\geq 4$ at the origin. In that
paper the authors implicitly conjectured that the base locus
of $|F_4|$ equals $C$, which was subsequently proved by G. Welters
\cite{welt}. Our proof follows from the fact that $|F_4|$ contains
all squares of quadrics in $|I(2)|$.

\bigskip

This paper leaves many questions unanswered (section 7), like e.g.
finding explicit equations of the quartics $F_W$, their syzygies,
the dimensions of $|F_3|$ and $|F_4|$. The techniques used here
also apply when replacing $|\omega|^*$ by Prym-canonical space
$|\omega \alpha|^*$, and generalizing rank-$2$ vector bundles to
symplectic bundles.

\bigskip

{\bf Acknowledgements:} Many results contained in this paper arose from
discussions with Bert van Geemen, whose influence on this work is
considerable. I would like to thank him for these enjoyable and
valuable discussions.

\section{Some constructions for rank-$2$ vector bundles with ca- nonical
determinant}

In this section we briefly recall some known results from \cite{bv}, 
\cite{vgi} and \cite{pp} on rank-$2$ vector bundles over $C$.

\subsection{Bundles $E$ with $\dim H^0(C,E) \geq 3$}

Let $W \subset H^0(C,\omega)$ be a $3$-plane. We denote by $[W] \in 
\Gr(3,H^0(\omega))$ the corresponding point in the Grassmannian and
by $\B \subset \Gr(3,H^0(\omega))$ the codimension $2$ subvariety
consisting of $[W]$ such that the net $\PP W \subset |\omega|$ 
has a base point. For $[W] \notin \B$ we consider (see \cite{vgi}
section 4) the rank-$2$ vector bundle $E_W$ defined by the exact sequence
\begin{equation} \label{esw}
0 \lra E^*_W \lra \cO_C \otimes W \map{ev} \omega \lra 0.
\end{equation}
Here $E^*_W$ denotes the dual bundle of $E_W$. We have $\det E_W =
\omega$ and $W^* \subset H^0(C,E_W)$. We denote by $\D$ the 
effective divisor in $|\cO_{\Gr}(g-2)|$ defined by the condition
$$ [W] \in \D \iff \dim H^0(C,E_W) \geq 4.$$
Moreover $\B \subset \D$ and $E_W$ is stable if $[W] \notin \D$.

\bigskip

Let $W^\perp \subset H^0(\omega)^* = H^1(\cO)$ denote the
annihilator of $W \subset H^0(\omega)$. We call the projective
subspace $\PP W^\perp \subset |\omega|^*$ the {\em vertex} and
denote by
$$ \pi: |\omega|^* \dashrightarrow \PP W^*, \qquad \pi:
C \rightarrow \Gamma \subset \PP W^*,$$
the projection with center $\PP W^\perp$. If $[W] \notin \B$,
then $C \cap \PP W^\perp = \emptyset$ and $\pi$ restricts to
a morphism $C \ra \PP W^*$. Its image is a plane curve $\Gamma$ 
of degree $2g-2$. We note that 
$E_W = \pi^* (T(-1))$, where $T$ is the tangent bundle of
$\PP W^* = \PP^2$.

\bigskip

Conversely any bundle $E$ with $\det E = \omega$ and $\dim
H^0(C,E) \geq 3$ is of the form $E_W$.

\subsection{Bundles $E$ with $\dim H^0(C,E) \geq 4$}

Following \cite{bv} (see also \cite{pp} section 5.2) we associate
to a bundle $E$ with $\dim H^0(C,E) = 4$ a rank $\leq 6$ quadric
$Q_E \in |I(2)|$, which is defined as the inverse image of the
Klein quadric under the dual $\mu^*$ of the exterior product map
$$ \mu^*: |\omega|^* \lra \PP (\Lambda^2 H^0(E)^*) \supset
\Gr(2, H^0(E)^*), \qquad Q_E:=(\mu^*)^{-1} \left( \Gr \right).$$
Composing with the previous construction, we obtain a rational
map
$$ \alpha : \D \dashrightarrow |I(2)|, \qquad \alpha([W]) = Q_{E_W}.$$
Moreover given a $Q \in |I(2)|$ with $\mathrm{rk} \ Q \leq 6$ and
$\mathrm{Sing} \ Q \cap C = \emptyset$, it is easily shown
that 
$$ \alpha^{-1}(Q) = \{ [W] \in \D \ | \ \PP W^\perp \subset Q \}.$$
If $\mathrm{rk} \ Q = 6$, then $\alpha^{-1}(Q)$ has two connected
components, which are isomorphic to $\PP^3$.

\begin{lem} \label{quaw}
We have $[W] \notin \D$ if and only if the linear map induced by restricting 
quadrics to the vertex $\PP W^\perp$ 
$$res: I(2) \lra H^0(\PP W^\perp,\cO(2))$$ 
is an isomorphism. 
\end{lem}

\begin{proof}
It is enough to observe that the two spaces have the 
same dimension and that a nonzero
element in $\ker res$ corresponds to a $Q \in |I(2)|$ with
$\rk Q \leq 6$.
\end{proof}

\subsection{Definition of the quartic $F_W$}

We will now define the main object of this paper. Given $[W]
\notin \B$, we consider the $2\theta$-divisor $D(E_W) \subset
JC$ (see e.g. \cite{bv},\cite{vgi},\cite{pp}), whose 
set-theoretical support equals
$$ D(E_W) = \{ \xi \in JC \ | \ \dim H^0(C, \xi \otimes E_W) > 0 \}.$$
Since $\mult_\cO D(E_W) \geq \dim H^0(C, E_W) \geq 3$ and since
any $2\theta$-divisor is symmetric, the first nonzero term of the Taylor
expansion of a local equation of $D(E_W)$ at the origin $\cO$ is a 
homogeneous polynomial $F_W$ of degree $4$. The hypersurface
in $|\omega|^* = \PP T_\cO JC$ associated to $F_W$ is also denoted by $F_W$. 
Here we restrict
attention to the case $\dim H^0(C,E_W) = 3 \  \text{or} \ 4$.  
We have 
$$ F_W := \mathrm{Cone}_\cO (D(E_W)) \subset |\omega|^*.$$

The study of the quartics $F_W$ for $[W] \in \Gr(3,H^0(\omega)) \setminus
\D$ is the main purpose of this paper. If $[W] \in \D$, the quartics
$F_W$ have already been described in \cite{pp} Proposition 5.12.

\begin{prop}
If $\dim H^0(C,E_W) = 4$, then $F_W$ is a double quadric
$$ F_W = Q^2_{E_W}.$$
\end{prop}

Since $|I(2)|$ is linearly spanned by rank $\leq 6$ quadrics (see \cite{pp}
section 5), we obtain the following fact, which will be used in section 6.

\begin{prop} \label{f4sq}
The linear subsystem $|F_4|$ contains all squares of quadrics in $|I(2)|$.
\end{prop}

Although we will not use that fact, we mention that the rational map
\eqref{fmap4} is given by a linear subsystem $\PP \Gamma \subset 
|\mathcal{J}_\B(g-1)|$, where $\mathcal{J}_\B$ is the ideal sheaf of
the subvariety $\B$. If $g=4$, the inclusion is an equality (see \cite{opp}
section 6). If $g>4$, a description of $\PP \Gamma$ is not known.

\section{Kempf's cohomological obstruction theory}

In this section we outline Kempf's deformation theory \cite{kempf1}
and apply it to the study of the tangent cones $F_W$ of the divisors
$D(E_W)$.

\subsection{Variation of cohomology}

Let $\EE$ be a vector bundle over the
product $C \times S$, where $S = \Spec(A)$ is an affine neighbourhood of the
origin of $JC$.
We restrict attention to the case
$$\EE = \pi_C^* E_W \otimes \mathcal{L},$$
for some $3$-plane $W$,
and recall that Kempf's deformation theory was applied \cite{kempf1},
\cite{kempf2}, \cite{ks} to the case $\EE = \pi_C^* M \otimes \mathcal{L}$,
for a line bundle $M$ over $C$. The line bundle $\mathcal{L}$ denotes
the restriction of a Poincar\'e line bundle over $C \times JC$ to the
neighbourhood $C \times S$. The fundamental idea to study the
variation of cohomology, i.e. the two upper-semicontinuous
functions on $S$
$$ s \mapsto h^0(C\times \{s\}, \EE \otimes_A \cc_s), \qquad
 s \mapsto h^1(C\times \{s\}, \EE \otimes_A \cc_s), $$
where $\cc_s = A/\m_s$ and  $\m_s$ is the maximal ideal of $s \in S$, is
based on the existence of an approximating homomorphism.

\begin{thm}[Grothendieck, \cite{kempf1} section 7] \label{gro}
Given a family $\EE$ of vector bundles over $C \times S$,
there exist two flat $A$-modules $F$ and $G$ of finite type and
an $A$-homomorphism $\alpha : F \ra G$ such that for all $A$-modules $M$,
we have isomorphisms
$$ H^0(C\times S , \EE \otimes_A M) \cong \ker(\alpha \otimes_A id_M),
\qquad
 H^1(C\times S , \EE \otimes_A M) \cong \coker(\alpha \otimes_A id_M).$$
\end{thm}

By considering a smaller neighbourhood of the origin, we may assume
the $A$-modules $F$ and $G$ to be locally free (Nakayama's lemma). Moreover
(\cite{kempf1} Lemma 10.2) by restricting further the neighbourhood,
we may find an approximating homomorphism $\alpha : F \ra G$ such that
$\alpha \otimes \cc_0 : F \otimes_A A/\m_0 \ra  G \otimes_A A/\m_0$ is the
zero homomorphism.

\bigskip

We apply this theorem to the family $\EE = \pi^*_C E_W \otimes \L$, for
$[W] \notin \D$.
Since by Riemann-Roch $\chi(\EE \otimes \cc_s) = \chi(E_W \otimes \L_s)
= 0$, $\forall s \in S$, and since $h^0(C,E_W) = 3$, the local equation
$f$ of the divisor
$$ D(E_W)_{|S} = \{ s \in S \: | \: h^0(C \times \{s\}, E_W \otimes \L_s) > 0
\} $$
is given at the origin $\cO$ by the determinant of
a $3 \times 3$ matrix of regular functions
$f_{ij}$ on $S$, with  $1 \leq i,j \leq 3$, which vanish at $\cO$, i.e.,
the $A$-modules $F$ and $G$ are free and of rank $3$. Hence
$$ f = \det (f_{ij}). $$

The linear part of the regular functions $f_{ij}$ is related to the
cup-product as follows (\cite{kempf1} Lemma 10.3 and Lemma 10.6): let
$\m = \m_0$ be the maximal ideal of the origin $\cO \in S$ and consider
the exact sequence of $A$-modules
$$ 0 \lra \m/\m^2 \lra A/\m^2 \lra A/\m \lra 0. $$
After tensoring with $\EE$ over $C \times S$ and taking cohomology,
we obtain a coboundary map
$$H^0(C,E_W) = H^0(C \times \{s \} , \EE \otimes_A A/\m )
\map{\delta}  H^1(C \times \{s \} , \EE \otimes_A \m/\m^2 ) =
H^1(C,E_W) \otimes \m/\m^2, $$
where $\m/\m^2$ is the Zariski cotangent space at $\cO$ to $JC$.
Note that we have a canonical isomorphism $(\m/\m^2)^* \cong H^1(\cO)$
and that a  tangent vector $b \in H^1(\cO)$ gives, by composing with the
linear form $l_b : \m/\m^2 \ra \cc$, a linear map $\delta_b : H^0(E_W) \ra
H^1(E_W)$. As in the line bundle case \cite{kempf1}, one proves

\begin{lem}
For any nonzero $b \in H^1(\cO) = T_{\cO} JC$, we have
\begin{enumerate}
\item
The linear map $\delta_b: H^0(E_W) \ra H^1(E_W)$ coincides with
the cup-product $(\cup b)$ with the class $b$, and is {\em
skew-symmetric} after identifying $H^1(E_W)$ with $H^0(E_W)^*$
(Serre duality).
\item
The coboundary map $\delta : H^0(E_W) \ra H^1(E_W) \otimes \m/\m^2$ is
described by a skew-symmetric $3 \times 3$ matrix $(x_{ij})$, with
$x_{ij} \in H^1(\cO)^*$. Moreover the linear form $x_{ij}$
coincides with the differential $(df_{ij})_0$ of $f_{ij}$ at the
origin $\cO$.
\end{enumerate}
\end{lem}
The coboundary map $\delta$ induces a linear map
$$ \Delta : H^1(\cO)  \lra  \Lambda^2 H^0(E_W)^*, \qquad
   b \lms  \delta_b,$$
which coincides with the dual of the multiplication map of global
sections of $E_W$. Moreover
$$\ker \Delta = W^\perp = \{ x_{12} = x_{13} = x_{23} = 0 \}. $$
Using a flat structure \cite{kempf2} we can write the power 
series expansion of the regular functions $f_{ij}$ around $\cO$
$$ f_{ij} = x_{ij} + q_{ij} + \cdots, $$
where $x_{ij}$ and $q_{ij}$ are linear and quadratic polynomials
respectively. We easily calculate the expansion of $f$: by
skew-symmetry its cubic term is zero, and its quartic
term equals
$$ F_W : q_{11} x^2_{23} + q_{22} x^2_{13} + q_{33} x^2_{12} + x_{12}x_{23}
(q_{13} + q_{31}) - x_{12} x_{23} (q_{12} + q_{21}) - x_{12} x_{13}
(q_{23} + q_{32}). $$
We straightforwardly deduce from this equation the following properties
of $F_W$.
\begin{prop} \label{singver}
\begin{enumerate}
\item The quartic $F_W$ is singular along the vertex $\PP W^\perp$.
\item For any $x \in W^\perp$, the cubic polar $P_x(F_W)$ is singular
along the vertex $\PP W^\perp$.
\end{enumerate}
\end{prop}

\subsection{Infinitesimal deformations of global sections of $E_W$}

We first recall some elementary facts on principal parts. Let
$V$ be an arbitrary vector bundle over $C$ and let $\Rat(V)$ be the
space of rational sections of $V$ and $p$ be a point of $C$.
The space of principal parts of $V$ at $p$ is the quotient
$$ \Prin_p(V) = \Rat(V)/\Rat_p(V), $$
where $\Rat_p(V)$ denotes the space of rational sections of $V$ which are
regular at $p$. Since a rational section of $V$ has only finitely
many poles, we have a natural mapping
\begin{equation}\label{prinpart}
 \mathrm{pp} : \Rat(V)  \lra  \Prin(V) := \bigoplus_{p \in C} \Prin_p(V),
\qquad
s  \lms  \left( s \  \mod \  \Rat_p(V) \right)_{p \in C}.
\end{equation}
Exactly as in the line bundle case (\cite{kempf1} Lemma 3.3), one proves
\begin{lem}
There are isomorphisms
$$ \ker \mathrm{pp}  \cong H^0(C, V), \qquad \coker \mathrm{pp}
\cong H^1(C, V). $$
\end{lem}

In the particular case $V = \cO$, we see that a tangent vector
$b \in H^1(\cO) = T_{\cO} JC$ can be represented
by a collection $\beta = \left( \beta_p \right)_{p \in I}$
of rational functions $\beta_p \in \Rat(\cO)$, where $p$
varies over a finite set of points $I \subset C$. We then define
$\mathrm{pp}(\beta) = \left( \omega_p \right)_{p \in I}
\in \Prin(\cO)$, where $\omega_p$ is the principal part of $\beta_p$ at $p$.
We denote by $[\beta] = b$ its cohomology class
in $H^1(\cO)$. Note that we can define powers of $\beta$ by
$\beta^k := \left( \beta_p^k \right)_{p \in I}$.

\bigskip
\noindent

For $i \geq 1$, let $D_i$ be the infinitesimal
scheme $\Spec(A_i)$, where $A_i$ is
the Artinian ring $\cc[\epsilon]/\epsilon^{i+1}$.
As explained in \cite{kempf2} section 2, a tangent vector $b \in H^1(\cO)$
determines a morphism
$$ \exp_{i,b} : D_i \lra JC,$$
with $\exp_{i,b}(x_0) = \cO$, where $x_0$ is the closed point of $D_i$. Let
$\LL_{i+1}(b)$ denote the pull-back of the Poincar\'e
sheaf $\L$ under the morphism $\exp_{i,b} \times id_C$.
Note that we have the following
exact sequences
\begin{equation} \label{esext1}
D_1 \times C : \qquad  0 \lra \epsilon \cO \lra \LL_2(b) \lra \cO \lra 0,
\end{equation}
\begin{equation} \label{esext2}
D_2 \times C : \qquad  0 \lra \epsilon^2 \cO \lra \LL_3(b) \lra
\LL_2(b)  \lra 0.
\end{equation}

The second arrows in each sequence correspond to the restriction to
the subschemes $\{x_0 \} \times C \subset D_1 \times C$ and
$D_1 \times C \subset D_2 \times C$ respectively.
As above we choose a representative $\beta$ of $b$. Following
\cite{kempf2} section 2, one shows that the space
of global sections $H^0(C \times D_i, \LL_{i+1}(b) \otimes E)$,
with $E = E_W$ and $[W] \notin \D$, is isomorphic to the
$A_i$-module
\begin{equation} \label{glosecllbeta}
V_i(\beta) = \{
f = f_0 + \cdots + f_i \epsilon^i \in \Rat(E) \otimes A_i \ \
\text{such that}
\
f \exp (\epsilon \beta) \ \text{is regular} \ \forall p \in C \}.
\end{equation}
An element $f \in V_i(\beta)$ is called an $i$-th order deformation
of the global section $f_0 \in H^0(E)$.
In the case $i=2$, the condition $f \in V_i(\beta)$ is equivalent
to the following three elements,
\begin{equation} \label{expdef}
f_0, \qquad f_1 + f_0 \beta, \qquad
f_2 + f_1 \beta + f_0 \frac{\beta^2}{2},
\end{equation}
being regular at all points $p \in C$ --- for $i=1$, we consider the
first two elements. Alternatively this means
that their classes in $\Prin(E)$ are zero. We note that,
given two representatives $\beta = \left( \beta_p \right)_{p \in I}$
and $\beta' = \left( \beta'_p \right)_{p \in I'}$ with $[\beta] =
[\beta']$, the two subspaces $V_i(\beta)$ and $V_i(\beta')$ of
$\Rat(E) \otimes A_i$ are different and that any rational
function $\varphi \in \Rat(\cO)$
satisfying $\mathrm{pp} (\varphi) = \mathrm{pp}
(\beta' - \beta)$ induces an isomorphism
$V_i(\beta) \cong V_i(\beta')$.

\bigskip
\noindent

We consider a class $b \in H^1(\cO) \setminus W^\perp$ and
a representative $\beta$ such that $[\beta] = b$. By taking
cohomology of \eqref{esext1} tensored with $E$, we observe that
a first order deformation of $f_0$, i.e., a global section
$f = f_0 + f_1 \epsilon \in  V_1(\beta) \cong
H^0(C \times D_1, \LL_2(b) \otimes E)$
always exists. Since $\rk (\cup b) = 2$, the global section
$f_0$ is uniquely determined up to a scalar
$$ f_0 \cdot \CC = \ker \left( \cup b : H^0(E) \lra H^1(E)  \right).$$
Moreover any two first order deformations of $f_0$ differ
by an element in $\epsilon H^0(E)$.

\bigskip
\noindent

We now state a criterion for a tangent vector
$b = [\beta]$ to lie on the quartic tangent cone $F_W$ in terms
of a second order deformation of $f_0 \in H^0(E)$.

\begin{lem} \label{criterion}
A cohomology class $b = [\beta]  \in H^1(\cO) \setminus W^\perp$ is contained
in the cone  over the quartic $F_W$ if and only if there exists a
global section
$$f = f_0 + f_1 \epsilon + f_2 \epsilon^2 \in V_2(\beta) \cong
H^0(C \times D_2 , \LL_3(b) \otimes E).$$
\end{lem}

\begin{proof}
The proof is similar to \cite{ks} Lemma 4. We work over the Artinian
ring $A_4$, i.e., $\epsilon^5 = 0$.
By Theorem \ref{gro} applied to the family $\LL_5(b) \otimes E$ over
$C \times D_4$, there exists an
approximating homomorphism of $A_4$-modules
\begin{equation} \label{approxhomo}
A_4^{\oplus 3} \map{\varphi} A_4^{\oplus 3},
\end{equation}
such that $\ker \varphi_{|D_2} \cong H^0(C \times D_2, \LL_3(b) \otimes E)$,
$\coker \varphi_{|D_2} \cong H^1(C \times D_2, \LL_3(b) \otimes E)$,
and $\varphi \otimes \CC_0 = 0$. We denote by $\varphi_{|D_2}$ the
homomorphism obtained from \eqref{approxhomo} by projecting to $A_2$.
Note that any $A_4$-module is free. The matrix $\varphi$ is
equivalent to a matrix
$$ M:= \left( \begin{array}{ccc}
\epsilon^u & 0 & 0 \\
0 & \epsilon^v & 0 \\
0 & 0 & \epsilon^w
\end{array}
\right).
$$
Since $\varphi \otimes \CC_0 = 0$, we have $u,v,w \geq 1$. Moreover
we can order the exponents so that $1 \leq u \leq v \leq w$. It follows
from the definition of $D(E_W)$ as a determinant divisor that the
pull-back of $D(E_W)$ by $\exp_4: D_4 \lra JC$ is given by the
equation (in $A_4$)
$$ \det M = \epsilon^{u+v+w}.$$
We immediately see that $b \in F_W$ if and only if $u+v+w \geq 5$.
Let us now restrict $\varphi$ to $D_1$,i.e., we project \eqref{approxhomo}
to $A_1$. Since we assume $b \notin W^\perp = \ker \Delta$, the
restriction $\varphi_{|D_1}$ is nonzero and by skew-symmetry of rank $2$,
i.e., $u=v=1$ and $w \geq 2$. Hence $b \in F_W$ if and only if
$w \geq 3$.

\bigskip

On the other hand the $A_2$-module $\ker \varphi_{|D_2} \cong
H^0(C \times D_2, \LL_3(b) \otimes E)$ has length $2+w$. Let $\mu$
be the multiplication by $\epsilon^2$ on this $A_2$-module. Then
by \eqref{glosecllbeta} the $A_2$-module $\ker \mu$ is
isomorphic to the $A_1$-module $H^0(C \times D_1, \LL_2(b) \otimes
E)$, which is of length $4$, provided $b \notin W^\perp$. Hence
we obtain that $w \geq 3$ if and only if there exists an $f \in
H^0(C \times D_2, \LL_3(b) \otimes E)$ such that $\mu(f) = \epsilon^2
f_0$. This proves the lemma.

\end{proof}

\section{Study of the quartic $F_W$}

In this section we prove geometric properties of the quartic $F_W$.

\subsection{Criteria for $b \in F_W$}

We now show that the criterion of Lemma \ref{criterion} simplifies
to a criterion involving only a first order deformation $f = f_0 +
f_1 \epsilon  \in V_1(\beta)$ of $f_0$. As above we assume $b \notin
W^\perp$.

\bigskip

First we observe that the rational differential form $f_1 \wedge f_0$
is independent of the choice of the representative
$\beta$, i.e., $f_1 \wedge f_0$
only depends on the cohomology class $b =[\beta]$: suppose we take
$\beta' = \left( \beta_p \cdot \varphi \right)_{p \in I}$,
where $\varphi \in \Rat(\omega)$. Then $f_0$ and $f_1$
transform into $f'_0 = f_0$ and $f'_1 = f_1 + \varphi f_0$,
from which it is clear that $f'_1 \wedge f'_0 = f_1 \wedge f_0$.

\bigskip

Secondly one easily sees that $f_0 = \pi(b)$ (section 2.1) and that,
under the canonical identification $\Lambda^2 W^* = \Lambda^2 H^0(E) =
W$, the $2$-plane $H^0(E) \wedge f_0$ coincides with the intersection
$V_b := H_b \cap W$, where $H_b$ denotes the hyperplane determined by $b \in
H^1(\cO)$.

\bigskip

It follows from these two remarks that, given $b$ and $W$, the
form $f_1 \wedge f_0$ is well-defined up to a regular differential
form in $V_b \subset W$.

\begin{prop} \label{simpl}
We have the following equivalence
$$ b \in F_W \qquad \iff \qquad f_1 \wedge f_0 \in H_b.$$
\end{prop}

\begin{proof}
Since $f_1 \wedge f_0$ does not depend on $\beta$, we may choose a
$\beta$ with simple poles at the points $p \in I$. By Lemma \ref{criterion}
and relation \eqref{expdef} we see that $b \in F_W$ if and only if
the cohomology class $[f_1 \beta + f_0 \frac{\beta^2}{2}]$ is zero
in $H^1(E) /  \im (\cup b)$ --- we recall that $f_1$ is defined up
to $H^0(E)$.

\bigskip

First we will prove that $[f_0 \frac{\beta^2}{2}] \in
\im(\cup b)$. The commutativity of the upper right triangle of the
diagram (see e.g. \cite{kempf1})
$$
\begin{array}{ccccccc}
 & & & & H^0(E) &  & \\
 & & & & & &  \\
 & & & & \downarrow \cdot \frac{\beta^2}{2} & \searrow \ \cup
[\frac{\beta^2}{2}] &  \\
 & & & & & & \\
H^0(E) & \lra & H^0(E(2I)) & \lra & E(2I)_{|2I} & \lra & H^1(E) \\
 & & & & & & \\
 & & \cap & & \cap & \nearrow & \\
 & & & & & & \\
 & & \Rat(E) & \map{\mathrm{pp}} & \Prin(E) & &
\end{array}
$$
implies that $[f_0 \frac{\beta^2}{2}] = f_0 \cup [\frac{\beta^2}{2}]$.
Moreover the skew-symmetric cup-product map
$\cup b$
$$\cup b = \wedge \overline{b}: \ H^0(E) = W^*
\lra H^1(E) = W = \Lambda^2 W^*$$
identifies with the exterior product $\wedge \overline{b}$, where
$\overline{b}= \pi(b) \in W^*$. It is clear that
$\im (\cup b) = \im (\wedge \overline{b}) =
\ker (\wedge \overline{b})$, where $\wedge \overline{b}$ also
denotes the linear form
\begin{equation} \label{wedgef0}
\wedge \overline{b}: \ \Lambda^2 W^* \lra \Lambda^3 W^* \cong \CC.
\end{equation}
As already observed, we have $f_0 = \overline{b}$. Denoting by
$c \in W^*$ the class $\pi([\frac{\beta^2}{2}])$, we see that the relation
$(f_0 \wedge c) \wedge \overline{b} =
\overline{b} \wedge c \wedge  \overline{b} = 0$ implies that
$f_0 \cup [\frac{\beta^2}{2}] \in \ker (\wedge \overline{b}) =
\im (\cup b)$.

\bigskip

Therefore the previous condition simplifies to
$[f_1 \beta] \in \im(\cup b)$. We next observe that the linear form
$\wedge \overline{b}$
on $H^1(E)$ \eqref{wedgef0} identifies with the exterior product map
$$ H^1(E) \map{\wedge f_0} H^1(\omega) \cong \CC.$$
Since we have a commutative diagram
$$
\begin{array}{rccccc}
f_1 \in & H^0(E(I)) & \map{\cdot \beta} & \Prin(E) & \lra & H^1(E) \\
 & & & & & \\
 &  &  & \downarrow \ \wedge f_0 & & \downarrow
 \ \wedge f_0 \\
 & & & & & \\
f_1 \wedge f_0 \in & H^0(\omega) & \map{\cdot \beta} & \Prin(\omega) &
\lra & H^1(\omega),
\end{array}
$$
and since $f_1 \wedge f_0 \in H^0(\omega) \subset \Rat(\omega)$,
we easily see that the condition $[f_1 \beta] \in \im(\cup b)$ is
equivalent to $f_1 \wedge f_0 \in H_b =
\ker \left( \cup b : H^0(\omega) \lra H^1(\omega) \right).$

\end{proof}

In the following proposition we give more details on the element
$f_1 \wedge f_0 \in H^0(\omega)$. We additionally assume that
$\pi(b) \notin \Gamma$, which implies that the
global section $f_0 \in H^0(E)$
does not vanish at any point and hence determines an exact sequence
\begin{equation} \label{extes}
 0 \lra \cO \map{f_0} E \map{\wedge f_0} \omega \lra 0.
\end{equation}
The coboundary map of the associated long exact sequence
\begin{equation} \label{extcl}
 \cdots \lra H^0(\omega) \map{\cup e} H^1(\cO) \lra \cdots
\end{equation}
is symmetric and coincides (e.g. \cite{kempf1} Corollary 6.8)
with cup-product $\cup e$ with
the extension class $e \in \PP H^1(\omega^{-1}) = |\omega^2|^*$. Moreover
$\cup e$ is the image of $e$ under the dual of the multiplication map
\begin{equation} \label{mapext}
 H^1(\omega^{-1}) = H^0(\omega^2)^* \hookrightarrow
\sym^2 H^0(\omega)^*, \qquad e \lms \cup e.
\end{equation}
We note that $\mathrm{corank}(\cup e) = 2$ and that $\ker (\cup e) =
V_b$. Hence $(f_1 \wedge f_0) \cup e$ is well-defined.

\begin{prop} \label{quad}
If $\pi(b) \notin \Gamma$, then $f_1 \wedge f_0 \notin \ker (\cup e)$ and
we have (up to a nonzero scalar)
$$(f_1 \wedge f_0) \cup e = b \in H^1(\cO).$$
\end{prop}

\begin{proof}
We keep the notation of the previous proof.  The condition
$f_1 \wedge f_0 \in V_b$ implies that $f_1$ is a regular section
and, by \eqref{expdef}, that $f_0$ vanishes at the support of $b$, i.e.,
$\pi(b) \in \Gamma$. As for the equality of the proposition, 
we introduce the rank-$2$ vector bundle $\hat{E}$ which 
is obtained from $E$ by
(positive) elementary transformations at the points $p \in I$ and
with respect to the line in $E_p$ spanned by the nonzero
vector $f_0(p)$. Then we have $E \subset \hat{E} \subset E(I)$ and
$\hat{E}$ fits into the exact sequence
$$ 0 \lra E \lra \hat{E} \lra \cO_I \lra 0. $$
Moreover $f_1 \in H^0(\hat{E})$, which follows from condition
\eqref{expdef}.  We also have the following exact sequences
$$
\begin{array}{ccccccccccc}
 0 & \lra & \cO(I) & \lra & \hat{E} & \map{\wedge f_0} & \omega & \lra & 0
& \qquad & (\hat{e}) \\
 & & & & & & & & & &  \\
 & & \cup & & \cup & & \Vert & & & &  \\
 & & & & & & & & & &  \\
 0 & \lra & \cO & \map{f_0} & E & \map{\wedge f_0} & \omega & \lra
& 0 & \qquad & (e),
\end{array}
$$
and the extension class $\hat{e} \in H^1(\omega^{-1}(D))$ is obtained
from $e$ by the canonical projection $H^1(\omega^{-1}) \ra
H^1(\omega^{-1}(I))$.  Taking the associated long exact sequences, we
obtain
$$
\begin{array}{cccccc}
f_1 \in & H^0(\hat{E}) & \map{\wedge f_0} & H^0(\omega) & \map{\cup \hat{e}} &
H^1(\cO(I)) \\
 & & & & &  \\
 & \cup &  & \Vert &  & \uparrow  \ \pi_I \\
 & & & & &  \\
 & H^0(E) & \map{\wedge f_0} & H^0(\omega) & \map{\cup e} & H^1(\cO),
\end{array}
$$
where the two squares commute. This means that
$$ \pi_I \left( (f_1 \wedge f_0) \cup e \right) = (f_1 \wedge f_0) \cup
\hat{e} = 0. $$
Since $f_1 \wedge f_0$ does not depend on $\beta$ (nor on $I$), the
latter relation holds for any $I$ with $I = \supp \beta$.  Hence,
denoting by $\langle I \rangle$ the linear span in $|\omega|^*$ of the
support $I$ of $\beta$, we obtain
$$ (f_1 \wedge f_0) \cup e \in \bigcap_{I = \supp \beta}  \ker \ \pi_I
= \bigcap_{b \in \langle I \rangle} \langle I \rangle = b.$$
\end{proof}

\subsection{Geometric properties of $F_W$}

\begin{prop} \label{geomprop}
For any $[W] \notin \D$ we have the following
\begin{enumerate}
\item The quartic $F_W$ contains the canonical curve $C$,i.e., 
$F_W \in |I(4)|$.

\item The quartic $F_W$ contains the secant line $\overline{pq}$,
with $p \not= q$,
if and only if $\overline{pq} \cap \PP W^\perp \not= \emptyset$ or
$\dim W \cap H^0(\omega(-2p-2q)) >0$.

\item Let $\Sigma$ be the set of points $p$ at which the tangent line
$\TT_p(C)$ intersects the vertex $\PP W^\perp$. Then $\Sigma$ is empty for
general $[W]$ and finite for any $[W]$. Moreover any point $p \in
C \setminus \Sigma$ is smooth on $F_W$ and the embedded tangent space
$\TT_p(F_W)$ is the linear span of $\TT_p(C)$ and $\PP W^\perp$.
\end{enumerate}
\end{prop}

\begin{proof}
All statements are easily deduced from Proposition \ref{simpl}. Given a
point $p \in C$ we denote by $\pp_p \in \Prin_p(\cO)$
the principal part supported at $p$ of a rational function with a
simple pole at $p$. Then the class
$[\pp_p] \in H^1(\cO)$ is proportional to $i_\omega(p) \in |\omega|^*
= \PP H^1(\cO)$ and the section $f_0$ vanishes at $p$. Hence
$f_0 \pp_p \in \Prin(E)$ is everywhere regular and we may choose
$f_1 =0$. This proves part 1. See also \cite{pp}.

\bigskip

As for part 2, we introduce $\beta_{\lambda,\mu} = \lambda \pp_p + \mu \pp_q
\in \Prin(\cO)$ for $\lambda, \mu \in \CC$ and denote by $s_p$ and
$s_q$ the global sections $\pi([\pp_p])$ and $\pi([\pp_q])$, which
vanish at $p$ and $q$ respectively. Then one checks that $f_0 =
\lambda s_p + \mu s_q \in \ker (\cup [\beta_{\lambda,\mu}])$ and
$\mathrm{pp}(f_1) = \lambda \mu (s_q \pp_p + s_p \pp_q) \in \Prin(E)$.
With this notation the condition of Proposition \ref{simpl}
transforms into
\begin{equation} \label{resfw}
0 = l_{\lambda, \mu}(f_0 \wedge f_1) = \lambda \mu (\lambda^2
\gamma_p + \mu^2 \gamma_q),
\end{equation}
where $l_{\lambda, \mu}$ is the linear form defined by
$[\beta_{\lambda,\mu}] \in H^1(\cO)$.
The scalars $\gamma_p$ and $\gamma_q$ are the
values of the section $s_p \wedge s_q \in W \cap
H^0(\omega(-p-q))$ at $p$ and $q$ respectively. We now conclude noting that
$s_p \wedge s_q = 0$ if and only if $\overline{pq} \cap \PP W^\perp
\not= \emptyset$.

\bigskip

As for part 3, we first observe that the assumption $\Sigma = C$ implies
that the restriction $\pi_{|C} : C \ra \PP W^*$ contracts $C$ to a point,
which is impossible. Next we consider the tangent vector $t_q$ at $p$
given by the direction $q$. By putting $\lambda = 1$ and $\mu = \epsilon$,
with $\epsilon^2 = 0$, into equation \eqref{resfw} we obtain that
$t_q \in \TT_p(F_W)$ if and only if $\epsilon \gamma_p = 0$, i.e.,
$\pi(q) \in \TT_{\pi(p)}(\Gamma)$. Hence $\TT_p(F_W) =
\pi^{-1}(\TT_{\pi(p)}(\Gamma))$,
which proves part 3.

\end{proof}

\subsection{The cubic polar $P_x(F_W)$}

Firstly we deduce from Propositions \ref{simpl} and \ref{quad} a
criterion for $b \in P_x(F_W)$, with $x \in W^\perp$. Let $H_x$ be
the hyperplane determined by $x \in H^1(\cO)$. As above we
assume $b \notin W^\perp$ and $\pi(b) \notin \Gamma$, i.e.,
the pencil $V = V_b$ is base-point-free.

\begin{prop} \label{cricub}
We have the following equivalence
$$ b \in P_x(F_W) \qquad \iff \qquad f_1 \wedge f_0 \in H_x.$$
\end{prop}

\begin{proof}
We recall from section 4.1 that $\cup e$ induces a symmetric
isomorphism $\cup e : (V^\perp)^* \map{\sim} V^\perp$ and we
denote by $Q^* \subset \PP (V^\perp)^*$ and $Q \subset \PP V^\perp$
the two associated smooth quadrics. Note that $Q$ and $Q^*$ are
dual to each other. Combining Propositions \ref{simpl}, \ref{quad}
and \ref{singver} (1) we see that the restriction of the quartic
$F_W$ to the linear subspace $\PP V^\perp \subset |\omega|^*$
splits into a sum of divisors
$$ \left(F_W \right)_{|\PP V^\perp} = 2 \PP W^\perp + Q.$$
We also observe that $Q$ only depends on $V$ (and on $W$) and not on $b$.
Taking the polar with respect to $x \in W^\perp$, we obtain
$$ \left(P_x(F_W) \right)_{|\PP V^\perp} = 2 \PP W^\perp + P_x(Q).$$
Finally we see that the condition $b \in P_x(Q)$ is equivalent
to $f_0 \wedge f_1 = (\cup e)^{-1} (b) \in H_x$.
\end{proof}

We easily deduce from this criterion some properties of $P_x(F_W)$.

\begin{prop} \label{cubcc}
The cubic $P_x(F_W)$ contains the canonical
curve $C$, i.e., $P_x(F_W) \in |I(3)|$.
\end{prop}

\begin{proof}
We first observe that the two closed conditions of Proposition
\ref{cricub} are equivalent outside $\pi^{-1}(\Gamma)$. Hence they
coincide as well on $\pi^{-1}(\Gamma)$ and we can drop the
assumption $\pi(b) \notin \Gamma$. Now, as in the proof of
Proposition \ref{geomprop}(1), we may choose $f_1 =0$.
\end{proof}

\begin{prop}
We have the following properties
$$ \bigcap_{x \in W^\perp} P_x(F_W) = S_W \cup \PP W^\perp \cup
\bigcup_{n \geq 2} \Lambda_n,$$
$$ F_W \cap S_W = C \cup \Lambda_1, \qquad \text{and} \qquad
\bigcup_{n \geq 0} \Lambda_n \subset F_W,$$
where $S_W$ is an irreducible surface. For $n \geq 0$, we denote by
$\Lambda_n$ the union of $(n+1)$-secant $\PP^n$'s to the
canonical curve $C$, which intersect the vertex $\PP
W^\perp$ along a $\PP^{n-1}$. If $W$ is general, then $\Lambda_n =
\emptyset$ for $n \geq 2$ and $\Lambda_1$ is the union of
$2(g-1)(g-3)$ secant lines.
\end{prop}

\begin{proof}
We consider $b$ in the intersection of all $P_x(F_W)$ and we first
suppose that $\pi(b) \notin \Gamma$. Then by Propositions
\ref{simpl} and \ref{cricub} we have
$$ f_0 \wedge f_1 \in \bigcap_{x \in W^\perp} H_x  = W.$$
Hence we obtain that $\PP V^\perp \cap \bigcap_{x \in W^\perp}
P_x(F_W)$ is reduced to the point $(\cup e)(W) \in \PP V^\perp$.
On the other hand a standard computation shows that $S_W$ is the
image of $\PP^2$ under the linear system of the adjoint curves of
$\Gamma$. Hence $S_W$ is irreducible.

\bigskip

If $\pi(b) \in \Gamma$, we denote by $p_1,\ldots,p_{n+1} \in C$
the points such that $\pi(p_i) = \pi(b)$. Then $f_0$ vanishes at
$p_1,\ldots,p_{n+1}$. Since $f_1 \wedge f_0$ does not depend on
the support of $b$, we can choose $\supp b$ such that $p_i \notin
\supp b$. Then $f_1$ is regular at $p_i$ and we deduce that $f_1
\wedge f_0 \in H^0(\omega(-\sum p_i)) \cap W = V_b.$ Now any
rational $f_1$ satisfying $f_1 \wedge f_0 \in V_b = \im (\wedge
f_0)$ is regular everywhere, which can only happen when $f_0$
vanishes at the support of $b$. By uniqueness we have $\supp b
\subset \{ p_1,\ldots,p_{n+1} \}$ and $b \in \Lambda_n$. Note that
$\Lambda_0 = C$. This proves the first equality.

\bigskip

If $b \in F_W \cap S_W$, we have $f_1 \wedge f_0 \in W \cap H_b =
V_b$ and we conclude as above. Note that $\Lambda_1$ is contained
in $S_W$ and is mapped by $\pi$ to the set of ordinary double
points of $\Gamma$.
\end{proof}

For any $[W] \in \Gr(3,H^0(\omega)) \setminus \D$ we introduce the subspace
of $I(3)$
$$ L_W = \{ R \in I(3) \  | \ R  \ \ \text{is singular along the vertex} \ \
\PP W^\perp \}. $$ Then Propositions \ref{cubcc} and
\ref{singver}(2)  imply that $P_x(F_W) \in L_W$. More precisely,
we have

\begin{prop}
The restriction of the polar map of the quartic $F_W$ to its vertex
$\PP W^\perp$
$$ \mathbf{P}: \ \  W^\perp \lra L_W, \qquad x \lms P_x(F_W),$$
is an isomorphism.
\end{prop}

\begin{proof}
First we show that $\dim L_W = g-3$. We choose a complementary
subspace $A$ to $W^\perp$,i.e. $H^0(\omega)^* = W^\perp \oplus A$,
and a set of coordinates $x_1,\ldots,x_{g-3}$ on $W^\perp$ and
$a_1,a_2,a_3$ on $A$. This enables us to expand a cubic $F \in
S^3 H^0(\omega)$
$$ F=F_3(x) + F_2(x)G_1(a) + F_1(x)G_2(a) + G_3(a), \qquad
F_i \in \CC[x_1,\ldots,x_{g-3}], \ G_i \in \CC[a_1,a_2,a_3],$$
with $\deg F_i = \deg G_i = i$. Let $\mathcal{S}_A$ denote the
subspace of cubics singular along $\PP A$,i.e. $G_2 = G_3 = 0$.
We consider the linear map
$$ \alpha : I(3) \lra \mathcal{S}_A, \qquad F \lms F_3(x) + F_2(x)
G_1(a).$$
Since by Lemma \ref{quaw} any monomial $x_ix_j \in H^0(\PP W^\perp, \cO(2))$
lifts to a quadric $Q_{ij} \in I(2)$, we observe  that the monomials
$x_ix_jx_k$ and $x_ix_ja_l$, which generate $\mathcal{S}_A$, also
lift e.g. to $Q_{ij}x_k$ and $Q_{ij}a_l$ in $I(3)$. Hence
$\alpha$ is surjective and $\dim L_W = \dim \ker \alpha$ is easily
calculated. One also checks that this computation does not depend on $A$.

\bigskip

In order to conclude, it will be enough to show that $\mathbf{P}$
is injective. Suppose that the contrary holds,i.e., there exists
a point $x \in W^\perp$ with $P_x(F_W) = 0$. Given any base-point-free
pencil $V\subset W$ and any $b \in V^\perp$, we obtain by Proposition
\ref{cricub} that $f_0 \wedge f_1 \in H_x$. Since $\cup e:
(V^\perp)^* \map{\sim} V^\perp$ is an isomorphism, we see that for
$b \notin (\cup e)^{-1}(H_x)$ the element $f_0 \wedge f_1$ must be
zero. This implies that $b \in \Lambda$ and since $b$ varies in
an open subset of $|\omega|^*$, we obtain $\Lambda = |\omega|^*$,
a contradiction.

\end{proof}

\subsection{The quadric bundle associated to $F_W$}

Let $\tilde{\PP}^{g-1}_W  \ra |\omega|^*$ denote the blowing-up
of $|\omega|^*$ along the vertex $\PP W^\perp \subset |\omega|^*$.
The rational projection $\pi : |\omega|^* \dashrightarrow \PP^2 =
\PP W^*$ resolves into a morphism $\tilde{\pi} :
\tilde{\PP}^{g-1}_W  \rightarrow \PP^2$. Since $F_W$ is
singular along $\PP W^\perp$ (Proposition \ref{singver} (2)), the
proper transform $\tilde{F}_W \subset \tilde{\PP}^{g-1}_W$ admits
a structure of a quadric bundle $\tilde{\pi} : \tilde{F}_W \ra
\PP^2$.

\bigskip

The contents of Propositions \ref{geomprop} and \ref{cubcc} can be
reformulated in a more geometrical way.

\begin{thm} \label{mainthm}
For any $[W] \in \Gr(3,H^0(\omega)) \setminus \D$, the quadric bundle
$\tilde{\pi}: \tilde{F}_W \ra \PP^2$ has the following properties
\begin{enumerate}
\item Its {\em Hessian} curve is $\Gamma \subset \PP^2$.

\item Its {\em Steinerian} curve is the
(proper transform of the) canonical curve $C \subset
|\omega|^*$.

\item The rational {\em Steinerian map} $\mathrm{St} : \Gamma 
\dashrightarrow C$, which associates to a singular quadric 
its singular point, coincides with the adjoint map $\mathrm{ad}$ 
of the plane curve $\Gamma$. Moreover the closure of the image 
$\mathrm{ad}(\PP^2)$ equals $S_W$.
\end{enumerate}
\end{thm}

\begin{rem}
We note that Theorem \ref{mainthm} is analogous to the main
result of \cite{ks} (replace $\PP^2$ with $\PP^1 \times \PP^1$).
In spite of this striking similarity and the relation between
the two parameter spaces $\sing$ and $\Gr(3,H^0(\omega))$ 
(see \cite{pp}), we were unable to find a common frame for both
constructions.
\end{rem}

\section{The cubic hypersurface $\Psi_V \subset \PP^{g-3}$ associated to a
base-point-free pencil $\PP V \subset |\omega|$}

In this section we show that the symmetric cup-product maps
$\cup e \in \sym^2 H^0(\omega)^*$ (see \eqref{extcl}) arise as
polar quadrics of a cubic hypersurface $\Psi_V$, which will
be used in the proof of Theorem \ref{mainthm2}.

\bigskip

Let $V$ denote a base-point-free pencil of $H^0(\omega)$. We consider the
exact sequence given by evaluation of sections of $V$
\begin{equation} \label{esv}
0 \lra \omega^{-1} \lra \cO_C \otimes V \map{ev} \omega \lra 0.
\end{equation}
Its extension class $v \in \mathrm{Ext}^1(\omega,\omega^{-1}) \cong 
H^1(\omega^{-2}) \cong H^0(\omega^3)^*$ corresponds to the 
hyperplane in $H^0(\omega^3)$, which is the image of the
multiplication map 
\begin{equation} \label{hyp}
\im \left(V \otimes H^0(\omega^2)  \lra H^0(\omega^3) \right).
\end{equation}
We consider the cubic form $\Psi_V$ defined by
$$ \Psi_V : \sym^3 H^0(\omega) \map{\mu} H^0(\omega^3) \map{\bar{v}} \CC,$$
where $\mu$ is the multiplication map and $\bar{v}$ the linear form
defined by the extension class $v$. It follows from the description
\eqref{hyp} that $\Psi_V$ factorizes through the quotient
$$ \Psi_V : \sym^3 \VV \lra \CC,$$
where $\VV := H^0(\omega)/V$. We also denote by $\Psi_V \subset
\PP \VV$ its associated cubic hypersurface.

\bigskip

A $3$-plane $W \supset V$ determines a nonzero vector $w$ in 
the quotient $\VV = H^0(\omega)/V$ and a general $w$ determines
an extension \eqref{extes} --- recall that $W^* \cong H^0(E)$. Hence
we obtain an injective linear map $\VV \hookrightarrow H^1(\omega^{-1}),
w \mapsto e$, which we compose with \eqref{mapext}
$$\Phi: \VV \hookrightarrow H^1(\omega^{-1}) = H^0(\omega^2)^*
\hookrightarrow \sym^2 H^0(\omega)^*, \qquad w \mapsto e \mapsto \cup e.$$
Since $V \subset \ker (\cup e)$, we note that $\im \Phi \subset \sym^2
\VV^*$.

\bigskip

We now can state the main result of this section.

\begin{prop}
The linear map $\Phi: \VV \ra \sym^2 \VV^*$ coincides with the polar map
of the cubic form $\Psi_V$, i.e.,
$$ \forall w \in \VV, \qquad \Phi(w) = P_w(\Psi_V).$$
\end{prop}

\begin{proof}
This is straightforwardly read from the diagram
obtained by relating the exact sequences \eqref{esv} and \eqref{esw}
via the inclusion $V \subset W$. We leave the details to the reader.
\end{proof}

We also observe that, by definition of the Hessian hypersurface (see
e.g. \cite{dk} section 3), we have an equality among degree 
$g-2$ hypersurfaces of $\PP \VV = \PP^{g-3}$
\begin{equation} \label{hesspsi}
 \Hess(\Psi_V) = \D \cap \PP \VV, 
\end{equation}
where we use the inclusion $\PP \VV \subset \Gr(3,H^0(\omega))$.

\begin{rem}
We recall (see \cite{dk} (5.2.1)) that the Hessian and Steinerian of a cubic
hypersurface coincide and that the Steinerian map is a rational involution
$i$. In the case of the cubic $\Psi_V$, the involution
$$ i: \Hess(\Psi_V) \dashrightarrow \Hess(\Psi_V)$$
corresponds to the involution of \cite{bv} Propositions 1.18 and
1.19, i.e., $\forall w \in
\D \cap \PP \VV$, the bundles $E_w$ and $E_{i(w)}$ are related
by the exact sequence
$$ 0 \lra E_{i(w)}^* \lra \cO_C \otimes H^0(E_w) \map{ev} E_w \lra 0.$$
Since we will not use that result, we leave its proof to the reader. 
\end{rem}

\section{Base loci of $|F_3|$ and $|F_4|$}

Let us denote by $|F_3| \subset |I(3)|$ and $|F_4| \subset |I(4)|$
the linear subsystems spanned by the image of the rational maps
$\F_3$  and $\F_4$ respectively. Then we have the following

\begin{thm} \label{mainthm2}
The base loci of $|F_3|$ and $|F_4|$ coincide with the canonical curve
$C \subset |\omega|^*$.
\end{thm}

\begin{proof}
Let $b \in \Bs |F_3|$ and let us suppose that $b \notin C$. We consider
a base-point-free pencil $V \subset H_b$. With the notation of
section 5, we introduce the rational map
$$  r_b : \PP \VV \dashrightarrow \PP \VV, \qquad w \mapsto
r_b(w) = w', \qquad \text{with} \ \tilde{\Psi}_V (w,w', \cdot) = b,$$
where $\tilde{\Psi}_V$ is the symmetric trilinear form of $\Psi_V$.
We note (Proposition \ref{quad}) that, for $w \notin \PP(H_b/V)$, the
element $r_b(w)$ is collinear with the nonzero
element $f_0 \wedge f_1  \ \mod \ V$ and that $r_b$ is defined
away from the hypersurface $\Hess(\Psi_V)$, which we assume to be
nonzero. Since $b \in \Bs|F_3|$ we obtain by Proposition \ref{cricub}
that
$$r_b(w) = \left( \bigcap_{x \in W^{\perp}} H_x \right) \ \text{mod}
\ V  = W \ \text{mod} \ V = w.$$
Hence $r_b$ is the identity map (away from $\Hess(\Psi_V)$). 
This implies that $\tilde{\Psi}_V(w,w,\cdot) = b$
for any $w \in \PP \VV$, hence $\Psi_V = x_0^3$, where $x_0$ is the
equation of the hyperplane $\PP(H_b/V) \subset \PP \VV$. This in
turn implies that $\Hess(\Psi_V) = 0$, i.e., $\PP \VV \subset \D$.
Since for a general $[W] \in \Gr(3,H^0(\omega))$ the pencil $V = W 
\cap H_b$ is base-point-free, we obtain that a general $[W]$ lies on the
divisor $\D$, which is a contradiction.

\bigskip

As for $|F_4|$, we recall that the fact $\Bs |F_4| = C$ follows from 
\cite{welt}. Alternatively, it can also be deduced by noticing 
(see Proposition \ref{f4sq}) that $\Bs |F_4| \subset \Bs |I(2)|$. Hence, if
$C$ is not trigonal nor a plane quintic, we are done. In the other
cases, the result can be deduced from Proposition \ref{geomprop} ---
we leave the details to the reader.
\end{proof}

\section{Open questions}

\subsection{Dimensions}

The projective dimensions of the linear systems $|F_3|$ and $|F_4|$ are
not known for general $g$. The known values of $\dim |F_4|$ for a general curve
$C$ are given as follows (see \cite{pp}).
$$
\begin{array}{|c|c|c|c|c|}
\hline
 g & 4 & 5 & 6 & 7  \\ 
\hline
 \dim |F_4|  & 4 & 15 & 40 & 88  \\
\hline 
\end{array}
$$
The examples of \cite{pp} section 6 show that $\dim |F_4|$ 
depends on the gonality of $C$. Moreover it can be shown that 
$|F_4| \not= |I(4)|$.

\subsection{Prym-canonical spaces and symplectic bundles}

The construction of the quartic hypersurfaces $F_W$ admit various
analogues and generalizations, which we briefly outline.

\bigskip

{\bf{(1)}}
Let $P_\alpha := \mathrm{Prym}(C_\alpha/C)$ denote the Prym
variety of the \'etale double cover $C_\alpha \ra C$ associated 
to the nonzero $2$-torsion
point $\alpha \in JC$. Given a general $3$-plane $Z \subset 
H^0(C, \omega \alpha)$, we associate the rank-$2$ vector bundle
$E_Z$ defined by 
$$ 0 \lra E_Z^* \lra \cO_C \otimes Z \map{ev} \omega \alpha \lra 0.$$
By \cite{ip} Proposition 4.1 we can associate to $E_Z$ the
divisor $\Delta(E_Z) \in |2\Xi|$, where $\Xi$ is a symmetric 
principal polarization on $P_\alpha$. Its projectivized tangent
cone at the origin $0 \in P_\alpha$ is a quartic hypersurface
$F_Z$ in the Prym-canonical space $\PP T_0P_\alpha \cong |\omega
\alpha|^*$. Kempf's obstruction theory equally applies to the
quartics $F_Z$. We note that $F_Z$ contains the Prym-canonical
curve $i_{\omega \alpha}(C) \subset |\omega \alpha|^*$.

\bigskip

{\bf{(2)}} Let $W$ be a vector space of dimension $2n+1$, for $n \geq 1$. 
We consider a {\em general} linear map 
$$ \Phi :  \Lambda^2 W^* \lra H^0(C,\omega).$$
By taking the $n$-th symmetric power $\sym^n \Phi$ and using the
canonical maps $\sym^n (\Lambda^2 W^*) \ra \Lambda^{2n} W^* \cong W$ and
$\sym^n H^0(\omega) \ra H^0(\omega^{\otimes n})$, we obtain a 
linear map
$$ \alpha: W \lra H^0(\omega^{\otimes n}),$$
which we assume to be injective. We then define the rank $2n$ vector
bundle $E_\Phi$ by
$$ 0 \lra E^*_\Phi \lra \cO_C \otimes W \map{ev} \omega^{\otimes n} \lra 0.$$
The bundle $E_\Phi$ carries an $\omega$-valued symplectic form and
the projectivized tangent cone at $\cO \in JC$
to the divisor $D(E_\Phi)$ is a hypersurface $F_\Phi$ in $|\omega|^*$ 
of degree $2n+2$. Moreover $F_\Phi \in |I(2n+2)|$.

\flushleft{Christian Pauly \\
Laboratoire J.-A. Dieudonn\'e \\
Universit\'e de Nice-Sophia-Antipolis \\
Parc Valrose \\
06108 Nice Cedex 2 \\
France \\
e-mail: pauly@math.unice.fr}

\end{document}